\documentclass{amsart}

\usepackage{amsthm}
\usepackage{amscd} 
\usepackage{amssymb}  
\usepackage{tikz}
\usepackage{hyperref}

\allowdisplaybreaks

\numberwithin{equation}{section} 

\newtheorem{Thm}{Theorem} 
\newtheorem{Prop}[Thm]{Proposition} 
\newtheorem{Lem}[Thm]{Lemma} 
 
\newtheorem{Fact}[Thm]{Fact}

\theoremstyle{remark} 
\newtheorem{Rem}[Thm]{Remark}

\theoremstyle{definition} 
\newtheorem{Def}[Thm]{Definition}

\newcommand\de{\partial}
\newcommand{\R}{\mathbb{R}}

\newcommand{\bit}{\begin{itemize}}
\newcommand{\eit}{\end{itemize}}

\author{Paolo Salvatore}
\title{The Fulton Mac Pherson operad and the W-construction}

\begin{document}
\maketitle
\setcounter{secnumdepth}{0}
In this short note we explain in detail the construction of a $O(n)$-equivariant isomorphism of topological operads $F_n \cong WF_n$ ,
where $F_n$ is the Fulton Mac Pherson operad and $W$ is the Boardman-Vogt construction \cite{BV}.  The non-equivariant version was announced
in \cite{Sal}.
An explicit isomorphism for $n=1$ was constructed in \cite{Barber}.
\medskip

Recall that $F_n(k)$ is a compactification of the space of ordered configuration space of $k$ points in $\R^n$ modulo translations and positive dilations
\cite{GJ}.
We consider the unbased version of our operad in the sense that $F_n(0)= \emptyset$. 
We recall that $F_n(k)$ is a smooth manifold with faces, i.e. a manifold with corners such that any codimension $l$ stratum is the transverse intersections of $l$ strata of codimension 1.  The symmetric group $\Sigma_k$ acts smoothly and freely on $F_n(k)$.
The strata of codimension $l$ are indexed by rooted trees with $k$ leaves labelled by $1,\dots,k$ and $l$ internal edges, such that any vertex has at least two incoming edges, and exactly one outgoing.  
The leaves are edges with no source, and the root is the only edge with no target. All other edges are internal. 
We call such trees {\em nested trees on $k$ leaves}. We call the number of incoming edges of a vertex $v$ its valence, and denote it by $|v|$.
The following property is well known and crucial.
\begin{Fact}\label{ft}
The stratum $F(T)$ of $F(k)$ indexed 
by a nested tree $T$ is  canonically diffeomorphic to the product $\prod_v F_n(|v|)$ over all  vertices $v$ of $T$.  
\end{Fact}
Observe that any nested tree $T$ on $k$ leaves defines an operad composition in a topological operad $O$ of the form
$\prod_v O(|v|) \to O(k)$, where the product runs over all vertices $v$ of $T$.  
Strictly speaking the composition is defined up to isomorphism, and requires a planar embedding of the tree. Alternatively we can consider an operad as a functor from 
the category of finite sets. A similar story holds for the $\circ_i$-operations that generate all operad compositions. They correspond to the compositions induced by
trees with a single internal edge. For details we refer to \cite{MSS}.
The diffeomorphism of  Fact \ref{ft} is induced by the operad composition  $$\theta: \prod_v F_n(|v|) \to F_n(k)$$  associated to the tree $T$, and identifies the source with its image. 

\medskip

 An element of $WF_n(k)$ instead is described by a nested tree on $k$ leaves $T$, together with an element $x_v \in F_n(|v|)$ for each vertex $v$ of $T$, and a ``length''
 parameter  $t_e \in [0,1]$ for each internal edge $e$ of $T$.  The description is unique up to the equivalence relation generated by the following move: 
if $l_e=0$ for some edge $e$, then we can compose operadically  the labels of its source and target, and collapse such edge to a single vertex labelled by the composition.
The operad composition $a \circ_i b \in WF_n$ along a tree with a single internal edge $e$ is the tree obtained by grafting together the $i$-th leaf of the labelled tree of $a$  with the root of the labelled tree of $b$, and declaring that the resulting new internal edge has length $1$.  

\begin{Thm} \label{main}
There is a $O(n)$-equivariant isomorphism of topological operads $$\beta:F_n \cong WF_n $$
\end{Thm}

The geometric idea of the proof is that $WF_n(k)$ can be seen as a fattening of the manifold with corners $F_n(k)$, since $WF_n(k)$ decomposes as union 
$$WF_n(k)=\bigcup_{T} F_n(T) \times [0,1]^{l(T)},$$ where $l(T)$ is the number of internal edges of $T$, and also the codimension of the correponding stratum $F_n(T)$.
This shows that $$WF_n(k) \cong F_n(k) \cup ([0,1]\times \de F_n(k))$$ but the right hand side is $\Sigma_k$-equivariantly diffeomorphic to $F_n(k)$ by the equivariant collar neighbourhood theorem for manifolds with corners.
\begin{Lem} (Equivariant collar neighbourhood theorem)  \label{collar}
Let $M$ be a compact manifold with faces on which a compact Lie group $G$ acts smoothly.
Then there is a $G$-equivariant collar of the boundary of $M$, i.e. a $G$-equivariant smooth embedding $c: \partial M \times [0,2]   \to M$
such that $c(2,x)=x$. 
\end{Lem}
\begin{proof}
By Theorem 4.1 in \cite{Mel} for each face $M_i$ there is a $G$-invariant smooth function $\rho_i$ on $M$ such that $\rho_i^{-1}(0)=M_i$, and a $G$-invariant vector field  $V_i$, such that 
$V_i \rho_j =\delta_{ij}$. If we set $g_i=2-\rho_i$ then the flow of the vector field $V_i$ gives an embedding $c_i:[0,2] \times M_i \hookrightarrow M$ such that $c_i(2,x)=x$ and $g_i(c_i(x,t))=t$. 
Furthermore $c_j(c_i(y,t_i),t_j)=c_i(c_j(y,t_j),t_i)$ of $y \in M_i \cap M_j$, and so 
we have embeddings 
$$e_I:\cap_{i \in I} M_i \times [0,2]^{|I|} \hookrightarrow M$$ If $x=e_I(y,(t_i)_i)$ is not in the image of another $e_J$ with $J \neq I$, then $c(x,t)=e_I(y,( (1/2+t/4) t_i)_i)$. This embedding is piecewise smooth and can 
easily be smoothened.

\end{proof}

We apply the lemma to our case $M=F_n(k), \, G=O(k) \times \Sigma_k$ and write $c_k=c$.
From now on we suppress the index $n$ from the notation and write $F=F_n$.

We build inductively on the arity $k$ the homeomorphism 
$\beta_k: F(k) \cong WF(k)$. 

\medskip

In arity $k=2$  $\beta_2:F(2) = WF(2)$ is the canonical identifcation. 
We recall that $F(2)$ is $\Sigma_2$-equivariantly homeomorphic to $S^{n-1}$ with the antipodal action.

\medskip

At the next stage $F(3)$ is a manifold with boundary equipped with a free action of $\Sigma_3$.
The boundary $\de F(3)$ is the union of three copies of $F(2) \times F(2)$, corresponding to the 3 nested trees on 3 leaves
with an internal edge. 
The space $WF(3)$ is obtained as $WF(3)= F(3) \cup_{\{0\} \times \de F(3)} ([0,1] \times \de F(3) )$.

\medskip

\begin{Def}
The map $\beta_3: F(3) \to  WF(3)$ sends 
\bit
\item $\beta_3(y)=y \in F(3) \subset WF(3)$ if $y \notin Im(c)$.
\item $\beta_3(c(t,x))=c(2t,x) \in F(3) \subset WF(3)$ for $t \in [0,1]$ and $x \in \de F(3)$ 
\item $\beta_3(c(t,x))=(t-1,x_1,x_2)_T \in \de F(3) \times [0,1] \subset WF(3)$  for  $t \in [1,2]$ and $x=x_1 \circ_T x_2  \in \de F(3)$. \label{quazz}
\eit
\end{Def}
In the last expression $(t-1,x_1,x_2)_T$ indicates the labelled tree $T$ with internal edge of length $t-1$, running from a valence 2 vertex labelled $x_2$ to a valence 2 vertex labelled $x_1$. 
It is easy to see that $\beta_3$ is a $\Sigma_3 \times O(n)$-equivariant homeomorphism. 
It also respects the operad composition, since for $x_1,x_2 \in F(2)$, the composition in $WF$  of 
$\beta_2(x_1)=x_1$ and $\beta_2(x_2)=x_2$ along a nested tree $T$ on $3$ leaves with two vertices and an internal edge is 
the labelled tree with vertices labelled $x_1,x_2$ and the internal edge of length $1$. But this 
is $$\beta_3(x_1 \circ_T x_2)=\beta(c(2,x_1 \circ_T x_2))=(1,x_1,x_2)_T$$

\medskip

We construct inductively $\beta_{k}$ for $k > 3$. We first extend the collar embedding $c_{k-1}$
to an embedding  
$$c'_{k-1}:[0,3]\times \de F(k-1) \to WF(k-1)$$  defined by
$$c'_{k-1}(t,x): =
\begin{cases}
\beta_{k-1}(c_{k-1}(t-1,x)) {\rm \quad for \;}2 \leq t \leq 3 \\
c_{k-1}(t,x){\rm \quad for\;  }0 \leq t \leq 2
\end{cases}$$

Now let us define $\beta_{k}:F(k) \to WF(k)$. 

\medskip
\bit
\item If $y \notin  \stackrel{\circ}{Im(c_k)}$ then $\beta_k(y)=y$ 
\item If $y=c(t,w)$ with $0 \leq t \leq 1$ then $\beta_k(c(t,w))=c(2t,w) \in F(k) \subset WF(k)$.

\item If $y=c(t,w)$ with $1 \leq t \leq 2$ and $w=x \circ_T \bar{x}$, then 
$\beta_k(c(t,w))$ is described by a labelled tree in $WF(k)$  that is obtained by grafting two labelled trees along $T$, 
a ``lower'' tree related to $x$ and an ``upper'' tree related to $\bar{x}$,
with the new internal edge of length $t-1$.  
\eit
There are three subcases for each tree, and so $3\cdot 3=9$ possible cases in total. We consider the lower tree:
\begin{enumerate} 

\item \label{zero} If $x \notin \stackrel{\circ}{Im(c)}$ then the lower tree is a single vertex with label $x$. 
\item \label{uno} If $x=c(s,z)$, with $0 \leq s \leq 1$ then the lower tree is a single vertex with label $c(st,z)$. 
\item \label{due} If $x=c(s,z)$ with $1 \leq s \leq 2$ then the lower tree is $c'(s+t-1,z)$. 
\end{enumerate}
A similar description holds for the upper tree 
(with the replacement $x \mapsto \bar{x}, s \mapsto \bar{s}, z \mapsto \bar{z}$).
We suppressed from the notation the index of $c$ that is the arity of $x$ (resp. of $\bar{x}$).

The induction process continues defining $c'_{k-1}$ and $\beta_k$ for all $k>3$. 

\begin{Prop} \label{cinque}
the map $\beta_k$ is well defined, $\Sigma_k \times O(n)$-equivariant, and continuous.
\end{Prop}
\begin{proof}
The function $\beta_k$ is defined as a piecewise continuous function on some closed sets 
and so we need to check that the definitions are compatible for $t=0, t=1, s=0, s=1, \bar{s}=0, \bar{s}=1$. 
Now the equality $c(0,w)=c(2\cdot 0,w)$ settles the case $t=0$. For $t=1$ observe that $c(2 \cdot 1,w)=w=x \circ_T \bar{x}$ is equivalent to the labelled 
tree obtained by grafting $x$ and $\bar{x}$ together, with a new internal edge of length $t-1=0$. Notice that $x$ and $\bar{x}$ do not change 
since in (\ref{uno})  $c(st,z)=c(s,z)$ and in (\ref{due}) $c'(s+t-1)=c'(s)=c(s)$.
Let us consider the lower tree.
For $s=0$ the element $x=c(0,z)$ is sent to $c(0\cdot t ,z)=x$ in both (\ref{zero}) and (\ref{uno}). 
For $s=1$ the element $x=c(1,z)$ is sent to $c'(1+t-1,z)=c(t,z)$ in both (\ref{uno}) and (\ref{due}).
A similar compatibility holds for the upper tree.
We also have to check that Definition \ref{quazz} does not depend on the operadic decomposition of $w$.
But by iterated applications of the definition it turns out that if $w$ is the operadic composition of elements $x_i$ along a tree $T$, then for $1 \leq t \leq 2$
$\beta_k(c(t,w))$ is the labelled tree obtained by grafting with edges of length $t-1$ the trees associated to $x_i$ by (\ref{zero}), (\ref{uno}), (\ref{due}),
(with $x,s,z$ replaced by appropriate $x_i,s_i,z_i$ ), and so the result does not depend on the order of the composition operations producing $w$.

The equivariance follows from the construction.

\end{proof}

\begin{Prop}\label{sei}
The map $\beta_k$ respects the operad composition.
\end{Prop}
\begin{proof}
For two arbitrary elements $x,\bar{x}$ of the operad $F$
we have that $$\beta_k(x \circ_T \bar{x})=\beta(c(2,x \circ_T \bar{x}))$$ is the labelled tree connecting a lower tree and an upper tree by an internal edge
of length $t-1=2-1=1$. The upper tree is $x$ if $x \notin Im(c)$; it is $c(2t,z)$ if $x=c(t,z)$; and it is $$c'(s+2-1,z)=c'(s+1,z)=\beta(c(s,z))=\beta(x)$$ otherwise.
In all cases it is $\beta(x)$. Similarly the lower tree is $\beta(x')$, therefore $$\beta_k(x \circ_T \bar{x})= \beta(x) \circ_T \beta(\bar{x})$$ where the latter 
composition takes place in $WF$.
\end{proof}
\begin{Prop}\label{sette}
The map $\beta_k$ is a homeomorphism.
\end{Prop}
\begin{proof}
We prove this by induction. It is clear that $\beta$ restricts to a homeomorphism from $(F(k)-Im(c_k)) \cup c_k([0,1] \times \de F(k))$ to $F(k) \subset WF(k)$ 
and we need only to verify that it restricts to a bijection from $c_k([1,2] \times \de F(k))$ to $\overline{WF(k)-F(k)}$.  
 We know that the proposition is true for $k=3$. If it is true for $k$  then
$$c'_k:[0,3] \times \de F(k) \to WF(k)$$ is an embedding. We prove simultaneously by induction that $$c'(\{1+t\} \times \de F(k)) \subset WF(k)$$ contains exactly the labelled trees
in $\overline{WF(k)-F(k)}$ with maximum edge length equal to $t-1$, for $1 \leq t \leq 2$. Namely by definition $c'(1+t,w)=\beta_k(c(t,w))$, when $w=x \circ_T \bar{x}$, is given by 
a labelled tree with an edge of length $t-1$, and an upper tree is either a vertex,
or a tree of type  
$c'(s+t-1,z)$ that by inductive hypothesis has a maximum edge  
length $\leq s+t-3 \leq t-1$, and a lower tree that behaves similarly. 
Now given a labelled tree in $x \in WF(k)$ we can decompose it by cutting it along all edges of maximum length $t-1$, obtaining some subtrees $T_i$, 
and then write it as $x=c'(1+t,w)$, where $w$ is the composition of appropriate indecomposable elements $x_i \in F(k_i)-\de F(k_i)$ such that 
the operations (\ref{zero}), (\ref{uno}), (\ref{due}) on $x_i$ produce the trees $T_i$. This decomposition exists and is unique by inductive hypothesis.

\end{proof}
Propositions \ref{cinque}, \ref{sei} and \ref{sette} together prove Theorem \ref{main}.

\medskip

\begin{Rem}
It is known that $F(k)$ is a piecewise algebraic (PA) manifold \cite{LV}, and it has a PA-action of $\Sigma_k \times O(n)$.
Together with Michael Ching we proved that a version of the Lemma \ref{collar} holds with $M=F(k)$, $c$ piecewise algebraic, and $G=\Sigma_k \times O(n)$ by constructing 
piecewise algebraic embeddings $c_i$ similarly as in the proof of the Lemma.
The same proof shows that Theorem \ref{main} holds in the PA-category, giving an isomorphism $F(k) \cong WF(k)$ of piecewise algebraic manifolds,
that is $\Sigma_k \times O(n)$-equivariant.
\end{Rem}


\begin{thebibliography}{98}
\bibitem{Barber}D.A. Barber, A Comparison of Models for the Fulton-Macpherson Operads, Ph.D. Thesis, Sheffield, 2017
\bibitem{BV} M. Boardman and R. Vogt, Homotopy Invariant Algebraic Structures on Topological Spaces, Lect.N.Math. 347, Springer, 1973  
\bibitem{GJ} E. Getzler and J.Jones, Operads, homotopy algebra, and iterated integrals for double loop spaces, arXiv: hep-th/9403055  
\bibitem{LV} P. Lambrechts and I. Volic, Formality of the little N-disks operad. Mem. Amer. Math. Soc. 230 (2014), no. 1079
\bibitem{MSS} M. Markl, S. Shnider and J. Stasheff, Operads in algebra, topology and physics, Math. Surveys and Mon. 96, A.M.S., 2002
\bibitem{Mel} P. Albin and R. Melrose, Resolutions of smooth groups actions. Spectral theory and geometric analysis, 1-26, Contemp. Math., 535, Amer. Math. Soc., Providence, RI, 2011
\bibitem{Sal} P. Salvatore, Configuration spaces with summable labels, in: Cohomological Methods in Homotopy Theory, Prog. in Math. 196 (2001) pp. 375-396,
Birkh\"auser, Berlin.

\end{thebibliography}
\end{document}